\documentclass[12pt,letterpaper, leqno]{article}
\usepackage{amsmath}
\usepackage{amsfonts}
\usepackage{amssymb}
\usepackage{amsthm}
\usepackage{fullpage}

\def\AAA{\mathbb{A}}
\def\FF{\mathbb{F}}
\def\PP{\mathbb{P}}
\def\ZZ{\mathbb{Z}}

\def\calO{\mathcal{O}}
\newtheorem{theorem}{Theorem}[section]
\newtheorem{cor}[theorem]{Corollary}
\newtheorem{lemma}[theorem]{Lemma}
\newtheorem{conj}[theorem]{Conjecture}
\numberwithin{equation}{theorem}

\DeclareMathOperator{\red}{red}

\author{Alina Bucur and Kiran S. Kedlaya}
\title{The probability that a complete intersection is smooth}
\date{October 6, 2012}
\begin{document}

\maketitle

\begin{abstract}
Given a smooth subscheme of a projective space over a finite field, we compute the probability that its intersection with a fixed number of hypersurface sections of large degree is smooth of the expected dimension. This generalizes the case of a single hypersurface, due to Poonen. We use this result to give a probabilistic model for the number of rational points of such a complete intersection. A somewhat surprising corollary is that the number of rational points on a random smooth intersection of two surfaces in projective 3-space is strictly less than the number of points on the projective line. 
\end{abstract}

\section{Introduction and results}

\let\thefootnote\relax\footnotetext{2010 Mathematics Subject Classification: 14G15, 11M38.}

One of the classical Bertini theorems states
that if a subscheme $X$ of a projective space over an infinite field is smooth,
then so is any sufficiently general hyperplane section of $X$. It follows that if $\dim(X) = m$,
then for any $k \in \{1, \dots, m\}$, 
the intersection of $X$ with $k$ sufficiently general hyperplanes is smooth of dimension $m-k$.

No meaningful analogue of this assertion exists over a finite field: the set of available hyperplanes is finite, 
so the set of ``sufficiently general'' hyperplanes may be empty. Katz proposed to address this issue by considering
not just hyperplanes but hypersurfaces. 
For a single hypersurface
of large degree, the probability that the intersection with $X$ is smooth of the correct dimension was computed by
Poonen \cite[Theorem~1.1]{poonen}.
\begin{theorem}[Poonen] \label{T:poonen1}
Let $X$ be a smooth quasiprojective subscheme of dimension $m \geq 1$ of the projective space $\PP^n$ over a finite field $\FF_q$. For $d$ a nonnegative integer, let $S_d$ be the set of homogeneous
polynomials of degree $d$ on $\PP^n$. Let
$\mathcal{P}_d$ be the set of $f \in S_d$ for which the hypersurface $H_f$
defined by $f$ intersects $X$ in a smooth scheme of dimension $m-1$. Then
\begin{equation} \label{eq:poonen}
\lim_{d \to \infty} \frac{\#\mathcal{P}_d}{\#S_d} = \zeta_X(m+1)^{-1}.
\end{equation}
\end{theorem}
In particular, there are infinitely many hypersurfaces whose intersection with $X$
is smooth of dimension $m-1$; this corollary was established independently by Gabber
\cite[Corollary~1.6]{gabber}. 

The equality \eqref{eq:poonen} is predicted by the following heuristic argument.
The zeta function $\zeta_X$ is defined as the Euler product
\[
\zeta_X(s) = \prod_{x \in X^\circ} (1 - q^{-s \deg(x)})^{-1},
\]
in which $X^\circ$ denotes the set of closed points of $X$, and $\deg(x)$ denotes the degree of $x$ over 
$\FF_q$. One shows easily that the product converges absolutely at $s=m+1$.
The right side of \eqref{eq:poonen} is thus equal to the product of $(1-q^{-(m+1)\deg(x)})$
over all $x \in X^{\circ}$; this factor computes the probability that (in some local coordinates)
the value and partial derivatives of $f$ do not all vanish at $x$.
This is precisely the probability that $H_f \cap X$ is smooth of dimension
$m-1$ at $x$, provided that we follow the convention of \cite{poonen} that
a scheme is considered to be smooth of any dimension at a point it does not contain.
With this convention, $H_f \cap X$ is smooth if and only if it is smooth at each $x \in X^{\circ}$,
so \eqref{eq:poonen} asserts that the local smoothness conditions behave asymptotically as if they were completely independent, even though they are in fact only independent in small batches (see Lemma~\ref{L:low degree}).
A related phenomenon is that the probability that a positive integer $N$ is squarefree is
$6/\pi^2 = 1/\zeta(2) = \prod_p (1-p^{-2})$, where $1-p^{-2}$ is the probability that $N$ is not divisible by
$p^2$; in fact,
one can formulate a conjectural common generalization of this statement and Theorem~\ref{T:poonen1}
\cite[\S 5]{poonen}.

The purpose of this paper is to generalize Poonen's theorem
to the case of complete intersections. Namely, for $\mathbf{d} = (d_1,\dots,d_k)$ a $k$-tuple of
positive integers, write $S_{\mathbf{d}}$ for the product $S_{d_1}\times \cdots \times S_{d_k}$.
For each $\mathbf{f} = (f_1,\dots,f_k) \in S_{\mathbf{d}}$, write $H_{\mathbf{f}}$ for $H_{f_1} \cap \cdots \cap H_{f_k}$.
We are now interested in the probability that for $\mathbf{f} \in S_{\mathbf{d}}$ chosen uniformly,
$H_{\mathbf{f}} \cap X$ is smooth of dimension $m-k$.
For $H_{\mathbf{f}} \cap X$ to be smooth of dimension $m-k$ at $x$, either $f_1,\dots,f_k$ must not all vanish at $x$,
or they must all vanish and have linearly independent gradients.
This event occurs with probability 
$1-q^{-k \deg(x)} + q^{-k \deg(x)} L(q^{\deg(x)}, m, k)$, where
\[
L(q, m, k) = \prod_{j=0}^{k-1} (1 - q^{-(m-j)})
\]
denotes the probability that $k$ randomly chosen vectors in $\FF_q^m$ are linearly independent.
One thus expects that in the limit
as $d_1,\dots,d_k \to \infty$, the probability that $H_{\mathbf{f}} \cap X$ is smooth tends to the product of these local probabilities.

Before stating this as a theorem, we follow Poonen \cite[Theorem~1.2]{poonen} by introducing the
option to modify finitely many local conditions. Identify the set $S_{\mathbf{d}}$ with the sections
over $\PP^n$
of the vector bundle $\calO_{\PP^n}(\mathbf{d}) = \oplus_{i=1}^k \calO_{\PP^n}(d_i)$.
We may impose modified local conditions at a finite set of closed points by forming a subscheme $Z$ of $X$
supported at those points, then specifying the possible images of $\mathbf{f}$ under the restriction
map $H^0(\PP^n, \calO_{\PP^n}(\mathbf{d})) \to H^0(Z, \calO_Z(\mathbf{d}))$. 

In terms of these notations, we have the following theorem. We make the error bound explicit for possible use
in some applications, as in \cite{bdfl}, and to clarify the extent to which the degrees must simultaneously tend to infinity
(the largest degree must grow subexponentially compared to the smallest degree).

\begin{theorem} \label{T:main}
Let $X$ be a quasiprojective subscheme of dimension $m \geq 0$
of the projective space $\PP^n$ over some finite field
$\FF_q$ of characteristic $p$.
Let $\overline{X}$ denote the Zariski closure of $X$ in $\PP^n$.
Let $Z$ be a finite subscheme of $X$ for which $U = X \setminus Z$ is smooth
of dimension $m$, and define $z = \dim_{\FF_q} H^0(Z, \calO_Z)$. 
Choose an integer $k \in \{1,\dots,m+1\}$,
a tuple $\mathbf{d} = (d_1,\dots,d_k)$ of positive integers with $z \leq d_1 \leq \cdots \leq d_k$,
and a subset $T$ of $H^0(Z, \calO_Z(\mathbf{d}))$. (Note that $\#H^0(Z, \calO_Z(\mathbf{d})) = q^z$.) Put
\[
\mathcal{P}_{\mathbf{d}} = \{\mathbf{f} \in S_{\mathbf{d}}: 
H_{\mathbf{f}} \cap U \mbox{ is smooth of dimension
$m-k$, and } \mathbf{f}|_Z \in T \}.
\]
Then
\begin{multline} \label{eq:big prob}
\frac{\#\mathcal{P}_{\mathbf{d}}}{\# S_{\mathbf{d}}} = 
\frac{\#T}{q^z}
\prod_{x \in U^{\circ}} \left(1 - q^{-k\deg(x)}  + q^{-k \deg(x)} L(q^{\deg(x)}, m, k) \right)\\
+ O((d_1 - z + 1)^{-(2k-1)/m} + d_k^{m} q^{-d_1/\max\{m+1,p\}}),
\end{multline}
where the implied constant is an increasing function of $n,m,k, \deg(\overline{X}),z$.
\end{theorem}

Note that for $k > m$, we only admit $\mathbf{f}$ when $H_{\mathbf{f}} \cap U$ is empty, and the
product in \eqref{eq:big prob} equals $\prod_{x \in U^{\circ}} ( 1 - q^{-k\deg(x)}) = 
\zeta_U(k)^{-1}$. Consequently, the probability that a single hypersurface section of $X$ is smooth is the
same as the probability that $m+1$ hypersurface sections of $X$ have empty intersection; this is analogous to the
fact that the probability that a single positive integer is squarefree is the same as the probability that
two positive integers are relatively prime. (One can again formulate a conjectural common generalization; see
Section~\ref{sec:further}.) By contrast, for $1 < k < m+1$, numerical evidence suggests that 
the product over $x \in U^{\circ}$ is not a rational number. If that is correct, then this product cannot
admit an interpretation as
an evaluation of a zeta function, or more generally as an evaluation
of the $L$-function of an \'etale sheaf of rank $k$.

The proof of Theorem~\ref{T:main} is an extension of Poonen's proof of Theorem~\ref{T:poonen1}. 
It is a sieving argument that separately treats the points of $X^{\circ}$ of low, medium, and high degree
(where these ranges are defined in terms of $d_1,\dots,d_k$).
For points in the low range (including the points of $Z$), one shows that the local conditions
are indeed independent. For a single point in the middle range, one similarly shows that singularities manifest
with the probability predicted by the local factor. One no longer has independence
of these local conditions, but they together contribute so little to the product that they can be controlled by
crude estimates. For points in the high range, there are too many points to control by such crude arguments.
Instead, we use Poonen's clever device of writing the $f_i$ so as to partially decouple
the low-order Taylor coefficients; one then bounds the effect of the points of high degree 
using B\'ezout's theorem. (This trick is the cause of the explicit appearance of $p$ in the error term,
as it relies on the
fact that the derivative of a $p$-th power vanishes in characteristic $p$.)

Note that Theorem~\ref{T:main} is not needed to deduce the existence of $k$ hypersurfaces
whose joint intersection with $X$ is smooth of dimension $m-k$, as this follows by induction
from Poonen's original theorem. Our intended application of Theorem~\ref{T:main}
is to compute the distribution of the
number of rational points on a random smooth complete intersection, via the following immediate
corollary of Theorem~\ref{T:main}.

\begin{cor} \label{C:bdfl}
Let $X$ be a smooth quasiprojective subscheme of dimension $m$ of
the projective space $\PP^n$ over a finite field
$\FF_q$. 
Let $\overline{X}$ denote the Zariski closure of $X$ in $\PP^n$.
Let $y_1,\dots,y_g,z_1,\dots,z_h$ be distinct $\FF_q$-rational points of $X$.
Choose an integer $k \in \{1,\dots,m\}$ and
a tuple $\mathbf{d} = (d_1,\dots,d_k)$ of positive integers with $g+h-1 \leq d_1 \leq \cdots \leq d_k$.
Then the probability that for $\mathbf{f} \in S_{\mathbf{d}}$,
$H_{\mathbf{f}} \cap X$ is smooth of dimension
$m-k$ and contains $y_1,\dots,y_g$ but not $z_1,\dots,z_h$ equals
\begin{multline}
\prod_{x \in X^\circ} \left( 1 - q^{-k\deg(x)} + q^{-k \deg(x)} L(q^{\deg(x)}, m, k) \right) \\
\left( \frac{q^{-k} L(q, m,k)}{1 - q^{-k} + q^{-k} L(q,m,k)} \right)^g \left( 
\frac{1-q^{-k}}{1 - q^{-k} + q^{-k} L(q,m,k)} \right)^h \\
+ O((d_1 - g-h+1)^{-(2k-1)/m}+
d_k^{m} q^{-d_1/\max\{m+1,p\}}),
\end{multline}
where the implied constant is an increasing function of $n,m,k, \deg(\overline{X})$.
\end{cor}

By Corollary~\ref{C:bdfl}, over a fixed field $\FF_q$, in the limit as all of the $d_i$ tend to infinity 
in such a way that $d_k^m q^{-d_1/(m+1)} \to 0$, the number of $\FF_q$-rational points
on a smooth complete intersection of hypersurfaces on $X$ of degrees $d_1,\dots,d_k$ is distributed
as the sum of $\#X(\FF_q)$ independent identically distributed Bernoulli random variables, each taking
the value 1 with probability $q^{-k} L(q, m,k)/(1 - q^{-k} + q^{-k} L(q,m,k))$ and 0 otherwise. One can also take the double limit as $d_1, \dots, d_k$ and $q$ all go to infinity, but only in a suitable range. The limiting distribution is, by the Central Limit Theorem, a standard Gaussian. We would like to remark that the limited range in which the CLT applies is probably an artifact of the shape of our error term, and not an intrinsic phenomenon.
In the case of plane curves (i.e., $X = \PP^n, n=2, k=1$), this was established by 
Bucur, David, Feigon, and Lal\'\i n \cite{bdfl} by essentially the method used here
(namely, by making explicit the implied error terms in \cite{poonen}). Analogous results had been previously
obtained for hyperelliptic curves by Kurlberg and Rudnick \cite{kurlberg-rudnick},
and for cyclic covers of $\PP^1$ of prime degree by Bucur, David, Feigon, and Lal\'\i n \cite{bdfl0};
in those cases, the error analysis is somewhat simpler (and the error bounds somewhat sharper) because
one can enforce smoothness simply by being careful about ramification.

Corollary~\ref{C:bdfl} has in turn the following consequence which we find mildly counterintuitive,
and which prompted the writing of this paper.
For a random smooth plane curve of degree $d$ in $\PP^2$, the average number of $\FF_q$-rational points
tends to 
\[
(q^2+q+1) \frac{q^{-1}(1-q^{-2})}{1 - q^{-3}} = q+1
\]
as $d$ tends to infinity; in other words, the average trace of Frobenius over smooth plane curves of degree $d$ in $\PP^2$ tends to 0.
(The same conclusion holds for a random hyperelliptic curve, by the trivial argument of pairing
each curve with its quadratic twist.)
On the other hand, for a random smooth intersection of two hypersurfaces of degrees $d_1, d_2$
in $\PP^3$, the average number of $\FF_q$-rational points
tends to 
\[
(q^3+q^2+q+1) \frac{q^{-2}(1-q^{-3})(1-q^{-2})}{1 - q^{-2} + q^{-2}(1-q^{-3})(1-q^{-2})} = 
q+1 - \frac{q^{-2}(1+q^{-1})}{1 + q^{-2} - q^{-5}}
\]
as $d_1,d_2$ tend to infinity appropriately. This limit is less than $q+1$, despite the fact that the range
$[0, q^3+q^2+q+1]$ for the number of points is much wider on the side greater than $q+1$!
The limit does tend to $q+1$ as $q$ grows, but for small $q$ the difference between the two
is large enough to be observed experimentally in computer simulations. (For instance, for $q=2$ the average is 
$37/13 < 2.847$.) The same sort of computation shows that the average number of points over $\FF_{q^2}$ of such a curve is 
\[q^2+q - q^{-1} + O\left(q^{-2}\right).\] This is consistent with Brock-Granville excess for $q$ large (see \cite{bg}).

More generally, for a random smooth intersection of $n-1$ hypersurfaces of degrees $d_1, \dots, d_{n-1}$ in $\PP^n$, in the limit as $d_1, \dots, d_{n-1} \to \infty$ appropriately, the
number of $\FF_q$-rational points is distributed
as the sum of $(q^{n+1}-1)/(q-1)$ i.i.d.~Bernoulli trials, each taking the value 1 with probability $q^{1-n} L(q, n,n-1)/(1 - q^{1-n} + q^{1-n} L(q,n,n-1))$ and 0 otherwise. 
 Hence the average number of $\FF_q$-points on such a smooth curve is equal to 
\[(q+1)- (q+1) (1 - q^{1-n}) \frac{1-  (1-q^{-n})\cdots (1-q^{-3})} {1 - q^{1-n} + q^{1-n}(1-q^{-n})\cdots (1-q^{-2})}<q+1,\] and it is of size $(q+1)\big(1+O(q^{-3})\big).$ For fixed $q$ and varying $n$, the average appears to decrease as $n$ increases
 (though we did not check this rigorously); it is easy to see that the limit as $n \to \infty$ is
\[
(q+1) \prod_{n=3}^{\infty} (1-q^{-n})
= q^{1 + 1/24} (1-q^{-1})^{-2} \eta(q^{-1}),
\]
where $\eta$ is the classical Dedekind eta function.

Even more generally, for a random smooth intersection of hypersurfaces of degrees $d_1,\dots,d_k$ in 
$\PP^n$, 
the average number of $\FF_q$-rational points tends to $q^{n-k} + q^{n-k-1} + \cdots +1$ if $k=1$, but to a limit strictly less than $q^{n-k} + q^{n-k-1} + \cdots +1$
if $k>1$. This can be seen as follows. One would get a limiting average of exactly $q^{n-k} + q^{n-k-1} + \cdots +1$
if the local condition for smoothness at a point $x$ were that the first-order Taylor approximations of $f_1,\dots,f_k$ had to be linearly independent. For $k=1$ this is the correct local condition,
but for $k>1$, this condition is too restrictive
when the sections do not all vanish at $x$. (One possible explanation is that for $k>1$, the intersection
of the hypersurfaces can be smooth without being geometrically integral. However, we suspect that this occurs
with probability 0 as the $d_i$ tend to infinity, and so does not account for the discrepancy.)

One can obtain other limiting distributions by considering families of different shapes.
For example, Wood \cite{wood} has shown that the average number of $\FF_q$-points on a
random trigonal curve over $\FF_q$ is slightly more than $q+1$, whereas the average number of $\FF_q$-points of a
random cyclic trigonal curve is exactly $q+1$ (as seen by averaging each curve with its cubic twists).
In a different direction, Kurlberg and Wigman \cite{kw} construct a sequence of surfaces in which (for $q$
fixed) the average number
of $\FF_q$-points of a random curve section becomes unbounded, and in fact the number of points
has a Gaussian distribution; this is arranged by choosing surfaces with many $\FF_q$-rational points, by taking
restrictions of scalars of curves which have many points over $\FF_{q^2}$ (e.g., classical or Drinfel'd modular curves). The example considered in the present paper is, to date, the only known family of curves that occurs naturally (starting with the projective space, and not some space with few points) for which the average number of points is less than $q+1.$

\subsection*{Acknowledgments}
Thanks to Nick Katz, Bjorn Poonen, and Ze\'ev Rudnick for helpful comments. We would also like to thank Melanie Matchett Wood for catching some errors in an earlier version of the paper.
Bucur was supported by NSF grant DMS-0652529 and
the Institute for Advanced Study (NSF grant DMS-0635607).
Kedlaya was supported by DARPA grant HR0011-09-1-0048,
NSF \mbox{CAREER} grant DMS-0545904,
MIT (NEC Fund, Cecil and Ida Green Career Development Professorship),
and the Institute for Advanced Study (NSF grant DMS-0635607, James D. Wolfensohn Fund).

\section{Proofs}

As noted earlier, the proof of Theorem~\ref{T:main} is a sieve over $X^{\circ}$ similar to 
Poonen's proof of Theorem~\ref{T:poonen1}, in which we separately analyze the contributions of points of low
degree (including the points of $Z$), medium degree, and high degree.
We analyze the points of low and medium degree using the following observation \cite[Lemma~2.1]{poonen}.
\begin{lemma} \label{L:low degree}
Let $Y$ be a finite closed subscheme of $\PP^n$ over $\FF_q$.
Then for any $d \geq \dim_{\FF_q} H^0(Y, \calO_Y)-1$, the restriction map
$\phi_d: H^0(\PP^n, \calO_{\PP^n}(d)) \to H^0(Y, \calO_Y(d))$ is surjective.
\end{lemma}

This lemma has the following corollaries, following \cite[Lemmas~2.2, 2.3]{poonen}.
\begin{cor} \label{C:low degree1}
For $r$ a positive integer,
let $U^{\circ}_{<r}$ be the subset of $U^{\circ}$ consisting of points of degree less
than $r$,
and define
\[
\mathcal{P}_{\mathbf{d}, r} = \{\mathbf{f} \in S_{\mathbf{d}}: H_{\mathbf{f}} \cap U \mbox{ is smooth of dimension $m-k$ at all
$x \in U^{\circ}_{<r}$, and $\mathbf{f}|_Z \in T$}\}.
\]
For
\begin{equation} \label{eq:bound on d1}
d_1 \geq (m+1) \sum_{x \in U^{\circ}_{<r}} \deg(x) + z - 1
\end{equation}
(where $z = \dim_{\FF_q} H^0(Z, \calO_Z) = \dim_{\FF_q} H^0(Z, \calO_Z(\mathbf{d}))$), we have
\begin{equation} \label{eq:big prob2}
\frac{\#\mathcal{P}_{\mathbf{d},r}}{\# S_{\mathbf{d}}} 
= 
\frac{\#T}{q^z}
\prod_{x \in U^{\circ}_{<r}} \left(1 - q^{-k\deg(x)} + q^{-k \deg(x)} L(q^{\deg(x)}, m, k) \right).
\end{equation}
\end{cor}
\begin{proof}
Apply Lemma~\ref{L:low degree} to the union of $Z$ with the subschemes defined by the square of the maximal
ideal at each point of $U^{\circ}_{<r}$.
\end{proof}

\begin{cor} \label{C:low degree2}
For $x \in U^{\circ}$ of degree $e \leq d_1/(m+1)$, the fraction of $\mathbf{f} \in S_{\mathbf{d}}$
for which $H_{\mathbf{f}} \cap U$ is smooth of dimension $m-k$ at $x$ equals
$1 - q^{-ke} + q^{-k e} L(q^{e}, m, k)$.
\end{cor}
\begin{proof}
Apply Lemma~\ref{L:low degree} to the subscheme defined by the square of the maximal ideal at $x$.
\end{proof}

The equation \eqref{eq:big prob2} accounts for the points of $U$ of low degree,
except that we need to control the difference between the products appearing on the right sides of  \eqref{eq:big prob}
and \eqref{eq:big prob2}. By doing so, we may also account for the points of medium degree.
\begin{lemma} \label{L:lang-weil}
For all positive integers $e$, $\#X(\FF_{q^e}) \leq 2^m \deg(\overline{X}) q^{me}$.
\end{lemma}
\begin{proof}
See \cite[Lemma~1]{lang-weil}.
\end{proof}

\begin{lemma} \label{L:medium degree}
Define
\begin{multline*}
\mathcal{Q}_{\mathbf{d}, r} = \{\mathbf{f} \in S_{\mathbf{d}}: H_{\mathbf{f}} \cap U \mbox{ is not smooth of dimension $m-k$}
\\ 
\mbox{ at some $x \in U^{\circ}$ with $r \leq \deg(x) \leq d_1/(m+1)$}\}.
\end{multline*}
Then $\#\mathcal{Q}_{\mathbf{d}, r}/\#S_{\mathbf{d}}$ is bounded above by $2^{m+1} \deg(\overline{X}) k q^{-r(2k-1)}$,
as is the difference between the products appearing on the right sides of \eqref{eq:big prob} and \eqref{eq:big prob2}.
\end{lemma}
\begin{proof}
The second assertion is obtained by first observing that
\[
L(q,m,k) \geq 1 - \sum_{j=0}^{k-1} q^{-m+j} \geq 1 - k q^{-m-k+1},
\]
and then calculating (using Lemma~\ref{L:lang-weil}) that
\begin{align*}
\sum_{x \in U^\circ \setminus U^{\circ}_{<r}}
q^{-k \deg(x)} - q^{-k \deg(x)} L(q^{\deg(x)}, m, k) &
\leq \sum_{e=r}^\infty \#U(\FF_{q^e}) k q^{-e(m+2k-1)} \\
&\leq \sum_{e=r}^\infty 2^m \deg(\overline{X}) q^{em} k q^{-e(m+2k-1)}\\
&= 2^m \deg(\overline{X})k \frac{q^{-r(2k-1)}}{1 - q^{-2k+1}} \\
&\leq 2^{m+1} \deg(\overline{X}) k q^{-r(2k-1)}.
\end{align*}
The first assertion follows from the same argument plus
Corollary~\ref{C:low degree2}.
\end{proof}

It remains to control the points of high degree. This reduces easily to the case $k=1$, which we treat
following Poonen \cite[Lemma~2.6]{poonen}.
\begin{lemma} \label{L:high degree}
Suppose that $U$ is contained in the affine space with coordinates $t_1,\dots,t_n$,
and $dt_1,\dots,dt_m$ freely generate $\Omega_{U/\FF_q}$. Then for any positive integers $d,e$, 
for $f \in S_d$ chosen uniformly at random,
the probability that there exists $x \in U^{\circ}$ for which $\deg(x) \geq e$ and $H_f \cap U$ is not smooth of dimension
$m-1$ at $x$
is at most \[(m+1) \deg(\overline X) \left((n-m)\deg(\overline{X})+ d-1\right)^m q^{-\min\{e, d/p\}}.\]
\end{lemma}
\begin{proof}
By dehomogenizing, we may identify $S_d$ with the set of polynomials of degree at most $d$ in $t_1,\dots,t_n$.
We may choose $f \in S_d$ uniformly at random by choosing $f_0 \in S_d$, 
$g_1,\dots,g_m \in S_{\lfloor (d-1)/p \rfloor}$, 
$h \in S_{\lfloor d/p \rfloor}$
uniformly at random, then putting
\[
f = f_0 + g_1^p t_1 + \cdots + g_m^p t_m + h^p.
\]
For $i=0,\dots,m$, define
\[
W_i = U \cap \left\{\frac{\partial f}{\partial t_1} = \cdots = \frac{\partial f}{\partial t_i} = 0\right\},
\]
defining the partial derivatives by writing $df$ in terms of $dt_1,\dots,dt_m$.
By choosing $n-m$ hypersurfaces of degree at most $\deg(\overline{X})$ which define linearly independent
relations among $dt_1,\dots,dt_n$ in $\Omega_{U/\FF_q}$ and then applying Cramer's rule, we see that
$\deg(\partial f/\partial t_i) \leq (n-m) \deg(\overline{X}) + (d-1)$.
Suppose $i \in \{0,\dots,m-1\}$, and 
fix a choice of $f_0,g_1, \dots,g_i$ for which $\dim(W_i) \leq m-i$. Let $V_1,\dots,V_{\ell}$ be the $(m-i)$-dimensional
irreducible components of $(W_i)_{\red}$. By B\'ezout's theorem, we have $\ell \leq \deg(\overline X)   \left((n-m) \deg(\overline{X}) + d-1\right)^i$. 
For each $j \in \{1,\dots,\ell\}$,
we have $\dim(V_j) \geq 1$, so we can choose $s \in \{1,\dots,m\}$ such that the projection
$x_s(V_j)$ has dimension 1. In particular, any nonzero polynomial in $x_s$ of degree at most
$\lfloor (d-1)/p \rfloor$ cannot vanish on $V_j$.
Since $\partial f/\partial t_{i+1} = \partial f_0/\partial t_{i+1} + g_{i+1}^p$,
the set of $g_{i+1}$ for which $\partial f/\partial t_{i+1}$ vanishes on $V_j$ is a coset of the
subspace of $S_{\lfloor (d-1)/p \rfloor}$ consisting of functions vanishing on $V_j$; we have just shown
that this subspace has codimension at least $\lfloor (d-1)/p \rfloor+1$. Hence the probability that $\partial f/\partial t_{i+1}$
vanishes on some component of $(W_i)_{\red}$ is at most $\ell q^{-\lfloor (d-1)/p \rfloor-1}
\leq  \deg(\overline X)   \left((n-m) \deg(\overline{X}) + d-1\right)^i q^{- d/p}$.

Now fix a choice of $f_0,g_1, \dots, g_m$ for which $W_m$ is finite. As in the previous paragraph,
$W_m^{\circ}$ consists of at most $\deg(\overline X) \left((n-m)\deg(\overline{X})+ d-1\right)^m$ points. If $x \in W_m^{\circ}$ has degree at least
$e$, then the set of $h \in S_{\lfloor d/p \rfloor}$ for which $f$ vanishes at $x$ is a coset of a subspace
of $S_{\lfloor d/p \rfloor}$ of codimension at least $\min\{e, \lfloor d/p \rfloor + 1\}$
\cite[Lemma~2.5]{poonen}.
Thus the probability that $H_f \cap U$ fails to be smooth at some point of degree at least $e$ is 
at most $\deg(\overline X) \left((n-m)\deg(\overline{X})+ d-1\right)^m q^{-\min\{e, d/p\}}$.
This yields the claim.
\end{proof}

\begin{cor} \label{C:high degree}
Define the set
\begin{multline*}
\mathcal{Q}_{\mathbf{d}}^{\mathrm{high}} = \{\mathbf{f} \in S_{\mathbf{d}}: H_{\mathbf{f}} \cap U \mbox{ is not smooth of dimension $m-k$} \\
\mbox{at some $x$ with $\deg(x) > d_1/(m+1)$}\}.
\end{multline*}
Then
\[
\frac{\#\mathcal{Q}_{\mathbf{d}}^{\mathrm{high}}}{\#S_{\mathbf{d}}} \leq 
k^m n^{2m} (m+1) \deg(\overline{X})^{m+1} d_k^{m} q^{-\min\{d_1/(m+1), d_1/p\}}.
\]
\end{cor}
\begin{proof}
We may reduce to the case $U \subseteq \AAA^n_k$ at the cost of multiplying by a factor of $n+1$ at the end.
Suppose that for some $i \in \{1,\dots,k\}$, $f_1,\dots,f_{i-1}$ have been chosen
so that $H_{f_1} \cap \cdots \cap H_{f_{i-1}} \cap U$ is smooth of dimension $m-i+1$ at each $x \in U^{\circ}$
with $\deg(x) \geq d_1/(m+1)$. 
Let $Y$ be the open subscheme obtained from $H_{f_1} \cap \cdots \cap H_{f_{i-1}} \cap U$ by removing all points
at which $H_{f_1} \cap \cdots \cap H_{f_{i-1}} \cap U$ fails to be smooth of dimension $m-i+1$; then
the Zariski closure $\overline{Y}$ of $Y$ has degree $\deg(\overline{X}) d_1 \cdots d_{k-1}$. 
We may cover $Y$ with at most $\binom{n}{m-i}$ open subsets, on each of which the module of differentials
is freely generated by some $(m-i)$-element subset of $dt_1,\dots,dt_n$.
We apply Lemma~\ref{L:high degree} to each of these sets, improving the bound slightly by noting 
that the relations among differentials used in Cramer's rule may be generated by first generating
relations for $U$, then adding the relations $df_1, \dots, df_{i-1}$.
With this, we see that for $f_i \in S_{d_i}$ chosen uniformly
at random, the probability that there exists $x \in U^{\circ}$ with $\deg(x) \geq d_1/(m+1)$
such that $H_{f_1} \cap \cdots \cap H_{f_i} \cap U$ is not smooth at $x$ is bounded above by
\begin{multline*}
\binom{n}{m-i} (m+1) \left(\deg(\overline X) d_1 \dots d_{i-1}\right) \left(n\deg(\overline{X}) + \sum_{j=1}^i (d_j-1)\right)^{m-i+1} 
q^{-\min\{d_1/(m+1), d_i/p\}} \\
\leq n^{2m} (m+1) k^m \deg(\overline{X})^{m+1} d_k^{m} q^{-\min\{d_1/(m+1), d_1/p\}}.
\end{multline*}
This yields the desired bound.
\end{proof}

It remains to complete the proof of Theorem~\ref{T:main}, by tuning the parameter $r$ appearing in
Corollary~\ref{C:low degree1}.
\begin{proof}[Proof of Theorem~\ref{T:main}]
Note that for $r$ a positive integer, by Lemma~\ref{L:lang-weil} we have
\begin{align*}
\sum_{x \in U^\circ_{<r}} \deg(x) &\leq \sum_{e=1}^{r-1} \#U(\FF_{q^e}) \\
&\leq \sum_{e=1}^{r-1} 2^m \deg(\overline{X}) q^{me} \\
&\leq 2^m \deg(\overline{X}) \frac{q^{m(r-1)}}{1 - q^{-m}} \\
&\leq 2^{m+1} \deg(\overline{X}) q^{m(r-1)}.
\end{align*}
Consequently, if $r$ is chosen so that
\[
d_1 - z + 1 \geq (m+1) 2^{m+1} \deg(\overline{X}) q^{m(r-1)},
\]
then \eqref{eq:bound on d1} will be satisfied.
We thus put
\[
r = 1 + \left\lfloor \frac{1}{m} \log_q \frac{d_1 - z + 1}{(m+1)2^{m+1} \deg(\overline{X})}
\right\rfloor;
\]
we may assume $r \geq 1$, as otherwise we may choose the implied constant in the big-O notation
so that the contribution of $(d_1 -z+1)^{-(2k-1)/m}$ bounds the error term.
By Lemma~\ref{L:medium degree} and Corollary~\ref{C:high degree},
\[
\frac{\#\mathcal{P}_{\mathbf{d}}}{\# S_{\mathbf{d}}} -
\frac{\#T}{q^z}
\prod_{x \in U^{\circ}} \left(1 - q^{-k\deg(x)}  + q^{-k \deg(x)} L(q^{\deg(x)}, m, k) \right)
\]
is bounded in absolute value by
\[
2^{m+2} \deg(\overline{X})k q^{-r(2k-1)}
+ k^m n^{2m} (m+1) \deg(\overline{X})^{m+1} d_k^{m} q^{-\min\{d_1/(m+1), d_1/p\}}.
\]
By virtue of our choice of $r$, this gives a bound of the desired form.
\end{proof}

To get Corollary \ref{C:bdfl} from Theorem \ref{T:main}, we follow the same procedure as in the proof of \cite[Proposition 1.5]{bdfl}. First, apply Theorem \ref{T:main} with $Z=\emptyset$ in order to compute 
\[\frac{ \# \{\mathbf{f} \in S_{\mathbf{d}}: 
H_{\mathbf{f}} \cap U \mbox{ is smooth of dimension
$m-k$}\}}{\#S_{\mathbf{d}}}.\]
Next, take $Z$ to be the disjoint union of one $\mathfrak m_{P}^2$-neighborhood $Z_P$ for each $P \in \{y_1, \ldots, y_g, z_1,\ldots z_h\}$ and 
identify $H^0(Z, \mathcal{O}_Z(\mathbf{d}))$ with $\prod_P H^0(Z_P, \mathcal{O}_{Z_P}(\mathbf{d}))$.
Take $T$ to be the set of tuples $(a_P)_P \in \prod_P H^0(Z_P, \mathcal{O}_{Z_P}(\mathbf{d}))$ with the following properties.
For $P = y_1,\dots,y_g$, we require that $a_P \in \mathfrak{m}_P H^0(Z_P, \mathcal{O}_{Z_P}(\mathbf{d}))$ (that is, all of the
sections vanish at $P$) and the individual factors of $a_P$ are linearly independent (that is, the intersection of 
the sections at $P$ is transversal).
For $P = z_1,\dots,z_h$, we require that $a_P \notin \mathfrak{m}_P H^0(Z_P, \mathcal{O}_{Z_P}(\mathbf{d}))$
 (that is, the sections do not all vanish at $P$).
Then apply Theorem~\ref{T:main} again to compute the probability that an element $\mathbf f \in S_{\mathbf d}$ gives an intersection $H_{\mathbf{f}} \cap X$ which is smooth of dimension $m-k$ passing through $y_1, \ldots, y_g$ and not through $z_1, \ldots, z_h.$
(Note that our choice of $T$ ensures that the intersection is smooth not only on $H_{\mathbf{f}} \cap U$ but also at each
point of $Z$.)
 Finally, taking the ratios of the two probabilities gives the desired result.  

\section{Distribution of number of rational points}

We mentioned at the end of the introduction that Theorem~\ref{T:main} can be used to study the distribution of the number of $\FF_q$-rational points on a random smooth complete intersection of $k$ hypersurfaces with $X$. Let us make this a bit more explicit, in the manner of \cite{bdfl}. Although there is some leeway to vary $X$ together with $q,d_1,\dots,d_k$, for simplicity we will instead fix a subscheme $\mathfrak{X}$ of the projective space over $\ZZ$, 
and always take $X$ to be the base change of $\mathfrak{X}$ to $\FF_q$. We assume further that each irreducible component
of the generic fibre of $\mathfrak{X}$ is geometrically irreducible; this guarantees that $\# X(\FF_q) \to +\infty$
as $q \to +\infty$. (One may instead fix a base finite field $\FF_{q_0}$, consider only those $q$ which are powers of
$q_0$, and take each $X$ to be the base change of a fixed subscheme $\mathfrak{X}$ of the projective space over $\FF_{q_0}$.
In this case, one should also assume that each irreducible component
of $\mathfrak{X}$ is geometrically irreducible.)


As stated before, the number of $\FF_q$-rational points
on a smooth complete intersection of hypersurfaces on $X$ of degrees $d_1,\dots,d_k$ is distributed, when the $d$'s get big enough,
as the sum of $\#X(\FF_q)$ independent identically distributed Bernoulli random variables, each taking
the value 1 with probability $q^{-k} L(q, m,k)/(1 - q^{-k} + q^{-k} L(q,m,k))$ and 0 otherwise.
The expected value of this sum of random variables is 
\[\#X(\mathbb F_q) \frac{q^{-k} L(q, m,k)}{1 - q^{-k} + q^{-k} L(q,m,k)},\]
and, since they are bounded real random variables, the normalized sum is distributed as a normal Gaussian when $q\rightarrow \infty$ by a suitable form of the central limit theorem \cite[Theorem~27.2]{billingsley}.

We would like to say the same about the number of points on the varieties in our family when
$d_1, \ldots , d_k$ and $q$ all go to infinity, which amounts to the computation of the moments. Denote
\[
N_r(q,{\bf d})=\frac{1}{\# S_{\bf d}^{\rm ns}} \sum_{{\bf f} \in S_{\bf d}^{\rm ns}}  \left(\frac{\#(H_{\bf f}\cap X)(\mathbb{F}_q)}{\sqrt{\# X(\FF_q) q^k} } \right)^r,
\]
where $S_{\bf d}^{\rm ns} = \{{\bf f}\in S_{\bf d}; H_{\bf f} \cap X \textrm{ is nonsingular}\}.$
We can write 
\begin{align*}
N_r(q,{\bf d}) 
=&\frac{q^{-rk/2} \#X(\FF_q)^{-r/2}}{\# S_{\bf d}^{\rm ns}} \sum_{{\bf f}\in S_{\bf d}^{\rm ns}}  \left(\sum_{y\in X(\mathbb{F}_q) }S_{\bf f}(y)\right)^r,
\end{align*}
where
\begin{equation} \label{eq:exp sum}
S_{\bf f}(y)=\frac{1}{q^k}\sum_{t \in \mathbb{F}_q}e\left( \frac{t f_1(y)}{q}\right ) \ldots \sum_{t\in \FF_q} e\left( \frac{tf_k(y)}{q}\right)=\begin{cases}1 & y \in H_{\bf f};\\
0 &\textrm{otherwise}. \end{cases}
\end{equation}
(Note that to give meaning to the evaluations $f_i(y)$, we must choose once for each $y \in X(\FF_q)$ a set of homogeneous coordinates defining $y$.)
Expanding the $r$-th power, 
\begin{align*}
N_r(q,{\bf d}) 
=&\frac{q^{-rk/2}  \#X(\FF_q)^{-r/2}  }{\# S_{\bf d}^{\rm ns}} \sum_{y_1, \ldots, y_r \in X(\mathbb{F}_q) } \sum_{{\bf f}\in S_{\bf d}^{\rm ns}}S_{\bf f}(y_1)\dots S_{\bf f}(y_r)\\
=& \frac{q^{-rk/2}  \#X(\FF_q)^{-r/2}}{\# S_{\bf d}^{\rm ns}} \sum_{g=1}^{\min\{r, \#X(\mathbb F_q)\} }a(g,r) \sum_{({\bf y}, {\bf b} ) \in P_{g,r}} \sum_{{\bf f}\in S_{\bf d}^{\rm ns}}S_{\bf f}(y_1)^{b_1}\dots S_{\bf f}(y_g)^{b_g},\\
\end{align*}
where 
\[P_{g,r}=\{((y_1,\ldots, y_g), (b_1, \ldots, b_g)); y_i \in X(\mathbb F_q) \textrm{ distinct, } b_i \in \mathbb Z_{>0}, b_1+\ldots +b_g=r \}, \]
while the $a(g,r)$ are certain combinatorial coefficients with the property that 
$$
\sum_{g=1}^r a(g,r)  \sum_{({\bf y}, {\bf b} ) \in P_{g,r}} 1 = \#X(\mathbb F_q)^r.
$$
For each $({\bf y}, {\bf b} ) \in P_{g,r}$, we have by \eqref{eq:exp sum}
\begin{equation}
\label{intsum}
\frac{1 }{\# S_{\bf d}^{\rm ns}} \sum_{{\bf f}\in S_{\bf d}^{\rm ns}}S_{\bf f}(y_1)^{b_1}\dots S_{\bf f}(y_g)^{b_g} = \frac{\#\{ {\bf f} \in S_{\bf d}^{\rm ns}; y_1, \ldots, y_g \in H_{\bf f} \} }{\# S_{\bf d}^{\rm ns}}.
\end{equation}
But Corollary \ref{C:bdfl}, together with the same trick employed in \cite[end of Section 2]{bdfl}, implies that the probability that $H_{\bf f} \cap X$ is smooth and passes through a specified set $y_1,\dots,y_g$ of distinct $\FF_q$-rational points of $X$ is equal to 
\begin{equation}\label{probmom}
\left( \frac{q^{-k} L(q, m,k)}{1 - q^{-k} + q^{-k} L(q,m,k)} \right)^g \left(1+
 O\left(q^{kg}\left(q^{-kr}(d_1 - g+1)^{-(2k-1)/m}+
d_k^{m} q^{-d_1/\max\{m+1,p\}}\right)\right)\right).
\end{equation}
Substituting \eqref{probmom} into \eqref{intsum}, then plugging the result into the formula for $N_r(q, {\bf d})$, we get that 
\begin{align*}
N_r(q, {\bf d})= & \frac{1}{q^{rk/2}\#X(\FF_q)^r}  \sum_{g=1}^{\min\{r, \#X(\mathbb F_q)\} }a(g,r) \sum_{({\bf y}, {\bf b} ) \in P_{g,r}} \left( \frac{q^{-k} L(q, m,k)}{1 - q^{-k} + q^{-k} L(q,m,k)} \right)^g \\
\times & \left(1+
 O\left(q^{k\min\{r, \#X(\mathbb F_q)\}}\left(q^{-kr}(d_1 - g+1)^{-(2k-1)/m}+
d_k^{m} q^{-d_1/\max\{m+1,p\}}\right)\right)\right).
\end{align*}
Since the main term equals the $r$-th moment of the sum of the aforementioned $\# X(\mathbb F_q)$ random variables divided by $q^{rk/2} \#X(\FF_q)^r$,  we conclude that the limiting distribution of the number of points on a smooth complete intersection $H_{\bf f} \cap X$ in $X$ over $\FF_q$ (normalized to have mean 0 and standard deviation 1) is a standard Gaussian, as long as $q$ and $d_1, \dots, d_k$ all go to infinity in such a way that the error term above approaches $0$. This condition is satisfied, for instance, if the $d_i$ are comparable in size (e.g., we cannot have $d_k$ exponentially bigger than $d_1$ because of the last term) and $d_1> q^{1+\epsilon}$ for any $\epsilon>0$. It is also satisfied if the characteristic $p$ of the base field stays bounded and, again, $d_k$ is of at most subexponential size compared to $d_1$.

\section{Further remarks} \label{sec:further}

We indicate some other directions in which one can probably generalize \cite{poonen} 
from hypersurfaces to complete intersections.

In the case of a single hypersurface section, one has an analogue
of Theorem~\ref{T:poonen1} allowing certain \emph{infinite} sets of modified local conditions
\cite[Theorem~1.3]{poonen}. The restriction on these local conditions is that at all but finitely
many points, they are no more stringent than the condition that $H_f$ must have smooth intersection
with each of a finite number of other varieties.
One expects to have a similar extension of Theorem~\ref{T:main}, in which one allows infinite sets of local
conditions provided that at all but finitely
many points, the conditions are no more stringent than the condition that $H_{\mathbf{f}}$ must have smooth intersection
with each of a finite number of other varieties.

One can also consider schemes of finite type over $\ZZ$, as in \cite[\S 5]{poonen}. We make the following
conjecture, generalizing the conditional theorem \cite[Theorem~5.1]{poonen}. One can again generalize by
allowing modification of 
finitely many local conditions; for simplicity, we omit these modifications from the following statement.

\begin{conj}
Let $X$ be a quasiprojective subscheme of $\PP^n_\ZZ$ which is regular of dimension $m \geq 0$.
For $d$ a positive integer, let $S_d$ be the set of homogeneous polynomials on $\PP^n_\ZZ$ of degree 
$d$, identified with $\ZZ^n$ using a $\ZZ$-basis of monomials.
For $\mathbf{d} = (d_1,\dots,d_k)$ a tuple of positive integers with $d_1 \leq \dots \leq d_k$,
put $S_{\mathbf{d}} = S_{d_1} \times \cdots \times S_{d_k}$. For $\mathbf{f} = (f_1,\dots,f_k) \in S_{\mathbf{d}}$,
put $H_{\mathbf{f}} = H_{f_1} \cap \cdots \cap H_{f_k}$, where $H_{f_i}$ denotes the hypersurface $f_i = 0$
on $\PP^n_\ZZ$. Put
\[
\mathcal{P}_{\mathbf{d}} = \{\mathbf{f} \in S_{\mathbf{d}}: 
H_{\mathbf{f}} \cap U \mbox{ is regular of dimension
$m-k$}\}.
\]
Then as $d_1,\dots,d_k \to +\infty$, the upper and lower densities of
$\mathcal{P}_{\mathbf{d}}$ in $S_{\mathbf{d}}$ both tend to
\[ 
\prod_{x \in X^{\circ}} \left(1 - q(x)^{-k}  + q(x)^{-k} L(q(x), m, k) \right),
\]
where $q(x)$ denotes the cardinality of the residue field of $x$.
\end{conj}

This conjecture includes the fact that the probability of a single integer being squarefree and the
probability of two integers being relative prime are both $6/\pi^2 = 1/\zeta(2)$. One may be able to prove
some special cases under the $abc$-conjecture, as is done in \cite{poonen}; the $abc$-conjecture is needed to 
compute the density of squarefree values of an integer polynomial. This is done in the univariate case by
Granville \cite{granville} and in the multivariate case by Poonen \cite{poonen-duke}.

\end{document}